\documentclass{article}
\usepackage{amssymb,amsmath}
\usepackage{amsfonts}
\usepackage{amsthm}
\usepackage{graphicx}

\newtheorem{Lem}{Lemma}
\newtheorem{The}{Theorem}
\newtheorem{Pro}{Proposition}
\newtheorem{Cor}{Corollary}

\theoremstyle{definition}

\newtheorem{Rem}{Remark}
\newtheorem{Exa}{Example}

\newtheorem{nota}{Notation}

\newcommand{\ds}{\displaystyle}

\newcommand{\for}{\qquad\mbox{for~}\,}

\def\RR{\mathbb{R}}
\def\ZZ{\mathbb{Z}}

\title{Periodic solutions for indefinite singular equations with applications to the weak case}

\author{Jos\'e Godoy\thanks{e-mail: jgodoy@ubiobio.cl}\\ {\small \itshape{Departamento de Matem\'atica,}}\\ {\small \itshape{Grupo de investigaci\'on en Sistemas Din\'amicos y Aplicaciones (GISDA),}}\\ {\small \itshape{Universidad del B\'io-B\'io, Casilla 5-C, Concepci\'on, Chile}} \and Manuel Zamora\thanks{e-mail: mzamora@ubiobio.cl}\\
{\small \itshape{Departamento de Matem\'atica,}}\\{\small \itshape{Grupo de Investigaci\'on en Sistemas Din\'amicos y Aplicaciones (GISDA)}}\\{\small \itshape{Universidad del B\'{\i}o-B\'{\i}o, Casilla 5-C, Concepci\'on, Chile}}}

\date{}

\begin{document}

\maketitle
\begin{abstract} 
As a consequence of the main result of this paper efficient conditions guaranteeing the existence of a $T-$periodic solution to the second order differential equation
\begin{equation*}
u''=\frac{h(t)}{u^{\lambda}}
\end{equation*}
are established. Here, $h\in L(\RR/T\ZZ)$ is a piecewise-constant sign-changing function where the non-linear term presents a weak singularity at 0 (i.e. $\lambda\in (0,1)$).  
\end{abstract}

\bigskip
{\small \noindent {\em MSC 2010 Classification} : 34C25, 34B18,
34B30

}

\section{Introduction and main results}

The study of the existence of $T-$periodic solutions to indefinite singular equations 
\begin{equation}\label{1}
u''=\frac{h(t)}{u^{\lambda}},
\end{equation}
where $\lambda>0$ and $h\in L(\RR/T\ZZ)$ has attracted the attention of many researchers during the last years. The motivation of such problems includes important applications in applied sciences: the stabilization of matter-wave breathers in Bose-Einstein condensates, the propagation of guided waves in optical fibers, or the electromagnetic trapping of a neutral atom near a charger wire (see \cite{BTans}).

Analytically, the above periodic boundary value problem, in spite of its simple looking structure, is considered a hard problem in the literature due to the lack of any a priori estimate over the set of (possible) periodic solutions, a condition used in general to apply one of the main tools of Nonlinear Functional Analysis: the Leray-Schauder degree theory.

Motivated initially by the pioneer paper of Bravo and Torres \cite{BTans}, the experts questioned whether this result could be extended to equation \eqref{1} under more general assumptions over $h$ and $\lambda$. However, this interesting equation still keeps some mystery as one can notice in \cite{Upre}, where the author shows that there is a sign-changing, $T-$periodic function $h\in C^{\infty}$ with mean value negative ($\overline{h}:=T^{-1}\int_0^T h(s)ds<0$) such that \eqref{1} with $\lambda=5/3$ has no $T-$periodic solution.

Therefore, the previous question should now be made more precise: what additional conditions on $\lambda$ and the weight function $h$ can be considered in order to ensure the existence of $T-$periodic solutions of \eqref{1}? (see, for e.g. \cite[Open problem 3.1]{HTjde}).

In \cite{Utmna}, the existence of a $T-$periodic solution was established when the weight function has only simple zeroes or is piecewise-constant with a finite number of values provided that $\lambda\geq 1$. In fact, a delicate relation between the order of the singularity and the multiplicities of the zeroes of the weight function was established in a recent paper \cite{HZpre2} when $\lambda\geq 1$, extending the previous result of Ure{\~n}a.

As it usually happens the situation can become even more difficult when weak singularities are considered (i.e., $\lambda\in (0,1)$). As far as we know this situation has not been considered until now.

Now we turn our attention to formulate the main result of this note which will be used to study equation \eqref{1} when $\lambda\in (0,1)$. With this proposal we denote by $\widetilde{L}(\RR/T\ZZ)$ the subspace of $L(\RR/T\ZZ)$ composed by the functions having mean value zero. Considering the equation
\begin{equation}\label{2}
u''=\frac{-\beta+\widetilde{h}(t)}{u^{\lambda}},
\end{equation}
where $\beta>0$ and $\widetilde{h}\in \widetilde{L}(\RR/T\ZZ)$, the main result of this paper is the following theorem.
\begin{The}\label{teo1}
For any non-trivial $\widetilde{h}\in \widetilde{L}(\RR/T\ZZ)$, there exists $\beta_*>0$ (depending on $\widetilde{h}$) such that \eqref{2} has at least one $T-$periodic solution for any $\beta\in (0,\beta_*]$.
\end{The}
Obviously the hypotheses of Theorem \ref{teo1} are necessary, because if we assume that $\beta$ is sufficiently small then the function defined by $h(t):=-\beta+\widetilde{h}(t)$ changes sign and $\overline{h}<0$, both conditions are needed in order that the periodic boundary value problem associated to \eqref{1} is solvable.

Theorem \ref{teo1} may be also applied to obtain information about the existence of a $T-$periodic solution of \eqref{1} according to the following general strategy: let $h\in L(\RR/T\ZZ)$ be a sign-changing function with $\overline{h}<0$, we define $\widetilde{h}(t):=h(t)-\overline{h}$ and apply Theorem \ref{teo1} in order to find $\beta_*(h)>0$ such that \eqref{2} has $T-$periodic solutions for any $\beta\in (0,\beta_*]$. To conclude that \eqref{1} has at least one $T-$periodic solution we need to assume that
\begin{equation}\label{3'}
\overline{h}\leq\beta_*(h)
\end{equation} 
(we choose $\beta:=-\overline{h}\in (0,\beta_*]$ in this case). However, since $\beta_*(h)$ depends on $h$, theoretically one may find a weight function $h\in L(\RR/T\ZZ)$ such that \eqref{3'} does not hold. Consequently, Theorem \ref{teo1} gives a general method to find $T-$periodic solutions provided that \eqref{3'} holds. This condition can replace the one established in the papers \cite{HZpre2, Utmna} when $\lambda\geq 1$.

Although this strategy can be widely used, in general to estimate the value of $\beta_*(h)$, in order to ensure \eqref{3'}, is a complicated task. However, the situation becomes more simple if $h$ is a piecewise-constant function. For instance, if
\begin{equation*}
h(t)=
\begin{cases}
h_1,&\quad {\mbox { if } } t\in [0,T/2),\\
-h_2,&\quad {\mbox { if } }t\in [T/2,T),
\end{cases}
\end{equation*}
where $h_1, h_2$ are positive constants and $\lambda\in (0,1)$ then an explicit corollary of Theorem \ref{teo1} is the following assertion.
\begin{Cor}\label{cor 1}
Assuming that 
\begin{equation}\label{hip}
0<\frac{h_2-h_1}{h_2+h_1}\leq\frac{\lambda(1-\lambda)}{4(1+\lambda)^2},
\end{equation}
then \eqref{1} has a $T-$periodic solution.
\end{Cor}
Corollary \ref{cor 1} has the following nice interpretation: if the forcing term $h$ is piecewise-constant with two values and $\lambda\in (0,1)$ then \eqref{1} has a $T-$periodic solution if one of the following items occur:
\begin{itemize}
\item $h$ is symmetrical (i.e., $h_1\sim h_2$);
\item $h$ has a large oscillation (i.e., $h_1+h_2$ is sufficiently large).
\end{itemize}
\begin{Exa}
Assume that $h$ is a piecewise-constant function defined as above. Then, the periodic problem
\begin{equation*}
u''=\frac{h(t)}{u^{1/2}},\qquad u(0)-u(T)=0=u'(0)-u'(T)
\end{equation*} 
has at least one solution provided that
\begin{equation*}
0<\frac{h_2-h_1}{h_2+h_1}\leq\frac{1}{36}.
\end{equation*}
\end{Exa}

The paper is structured as follows. In Section \ref{sec2} a change of variables is done in order to reduce our problem to an equivalent one which can be formulated as a manageable fixed point problem. The key result to prove Theorem \ref{teo1} is presented and proved in Section \ref{sec3} (see Lemma \ref{lema 2}) which is combined with some a priori bounds obtained in Section \ref{sec4} in order to prove Theorem \ref{teo1} and Corollary \ref{cor 1} in Section \ref{sec5}. Finally, in Section \ref{sec6} we extend the results to more general piecewise-constant forcing terms, and Section \ref{sec7} is devoted to prove two auxiliary lemmas used to optimize the value of the parameter $\beta_*>0$ given by Theorem \ref{teo1} under the hypotheses of Corollary \ref{cor 1}. For other results concerning periodic or Neumann solutions of nonlinear problems with a weight function having an indefinite sign see for instance \cite{BCNnodea, BGjdde, BZjde, BosZan, BFZpre, FGUjfa}.

For convenience, we finish this introduction with a list of notation which is used throughout the paper:
\begin{nota}\label{nota1}
Given $\widetilde{h}\in \widetilde{L}(\RR/T\ZZ)$ a non-trivial function, we denote
\begin{equation*}
\sigma(t):=\int_0^t\widetilde{h}(s)ds,\qquad \widetilde{H}:=\int_0^T\widetilde{h}^+(s)ds=\int_0^T\widetilde{h}^-(s)ds,
\end{equation*}
denoting by $b^+:=\max\lbrace 0,b\rbrace$, $b^-:=\max\lbrace 0,-b\rbrace$ for any real number $b$. We define
\begin{equation*}
\alpha:=\frac{1}{T}\int_0^T(\sigma(s)-\overline{\sigma})^2ds
\end{equation*}
and we fix $a>0$ such that $\sqrt{\lambda\alpha}<a$. Observe that the non-triviality of $\widetilde{h}$ implies that $\alpha>0$.

For any vector $\xi=(\beta,x,y)\in\RR^3$ and $f\in C([0,T];\RR)$, if appropriate, we shall denote
\begin{gather*}
A_{\xi}[f](t):=a+x+\beta\int_0^t f(s)ds,\for t\in [0,T]\\
m_f:=\min_{t\in [0,T]}f(t),\qquad M_{f}:=\max_{t\in [0,T]}f(t),\\
 B_{\xi}[f](t):=(a+x)y+\int_0^t\left[\left(\frac{2\lambda}{\lambda +1}\right)\beta f^2(s)+\left(\frac{1+\lambda}{2} \right)(-\beta+\widetilde{h}(s))\right]ds.
\end{gather*}
\end{nota}

\section{From a singular equation to an equivalent non-singular equation}\label{sec2}
Throughout this paper the family of equations
\begin{equation}\label{2.1}
u''=\frac{-\beta^2+\beta\widetilde{h}(t)}{u^{\lambda}},
\end{equation}
depending on a positive parameter $\beta>0$ and $\widetilde{h}\in\widetilde{L}(\RR/T\ZZ)$ a sign-changing function (i.e., $\widetilde{h}$ is non-trivial) will play a fundamental role. Generally speaking, it is a hard problem to study the existence of $T-$periodic solutions of \eqref{2} in a direct way. Hence, in this section we shall focus our attention to reduce the problem to an equivalent one where a deeper treatment can be done.

Observe that the change of variables
\begin{equation}\label{2.1.1}
\rho(t):=u^{\frac{\lambda+1}{2}}(t)\for t\in [0,T]
\end{equation}
allows us to pass to the equation
\begin{equation}\label{2.2}
\frac{d}{dt}\left(\rho'\rho\right)-\left(\frac{2\lambda}{1+\lambda}\right)\rho'^2(t)=\left(\frac{1+\lambda}{2}\right)(-\beta^2+\beta\widetilde{h}(t)).
\end{equation}
Notice that \eqref{2.2} takes the advantage of being a regular equation (it has no singularity) and therefore it should be easier to find new ideas or take a different approach to solve the problem. Indeed, the main result of this section states an equivalent formulation of the problem of finding positive $T-$periodic solutions of \eqref{2.2} (also \eqref{2.1}) noticing that it suffices to find points $(\beta,x,y,f)\in \RR^3\times C([0,T];\RR)$ satisfying certain properties. Under this formulation it will be more appropriate to apply the Leray-Schauder degree theory. 

\begin{Pro}\label{prop 2.1}
Assume that there exists $\xi=(\beta,x,y)\in\RR^+\times\RR^2$  and $f\in C([0,T],\RR)$  such that
\begin{gather}
\label{2.3} A_\xi[f](t)>0,\qquad A_\xi[f](t)f(t)=B_\xi[f](t)\for t\in[0,T],\\
\label{2.4} \frac{1}{T}\int_0^T f^2(s)ds=\frac{(1+\lambda)^2}{4\lambda},\qquad \int_0^T f(s)ds=0.
\end{gather}
Then $\rho(t):=A_{\xi}[f](t)$ is a positive $T-$periodic solution of \eqref{2.2}.
\end{Pro}
\begin{proof}
It follows immediately from the equations of \eqref{2.3} that $\rho(t)$ is positive and satisfies the differential equation. We only need to show that 
\begin{equation*}
\rho(0)=\rho(T),\qquad \rho'(0)=\rho'(T).
\end{equation*}
The first equality is obtained from the second equation of \eqref{2.4}. Notice that the first equation of \eqref{2.4} guarantees that
\begin{equation*}
B_{\xi}[f](0)=B_{\xi}[f](T),
\end{equation*}
which combined with
\begin{equation*}
f(t)=\frac{B_{\xi}[f](t)}{A_{\xi}[f](t)}
\end{equation*}
gives us that $f(0)=f(T)$. Since $\rho'(t)=\beta f(t)$, we obtain the second equality.
\end{proof}

In what follows we will work to find a connected set of solutions $(\beta,x,y,f)\in\RR^3\times C([0,T];\RR)$ to problem \eqref{2.3}-\eqref{2.4}. Of course, it will be convenient first to establish a functional analysis framework for our problem. 

\section{Continuation of solutions with nonzero degree}\label{sec3}

With the proposal of finding solutions $(\beta,x,y,f)\in\RR^3\times C([0,T];\RR)$ to the problem \eqref{2.3}-\eqref{2.4}, we now rewrite such problem as a fixed point equation depending on a parameter. In order to do this, we shall work on the Banach space $X:=\RR^2\times C([0,T];\RR)$ which is endowed with the usual norm $\|(x,y,f)\|:=|x|+|y|+\|f\|_{\infty}$. Problem \eqref{2.3}-\eqref{2.4} now becomes a fixed point equation on $X$ depending on the parameter $\beta:$
\begin{equation}\label{n1}
(x,y,f)=F[\beta,x,y,f],
\end{equation}
where
\begin{equation*}
F[\beta,x,y,f]:=\left(x-\int_0^T f(s)ds,y-\frac{1}{T}\int_0^T f^2(s)ds+\frac{(1+\lambda)^2}{4\lambda} ,\frac{B_{(\beta,x,y)}[f]}{A_{(\beta,x,y)}[f]}\right).
\end{equation*}
Solutions to \eqref{2.3}-\eqref{2.4} make sense on the set
\begin{equation*}
\Lambda:=\left\lbrace(\beta,x,y,f):A_{(\beta,x,y)}[f](t)>0,\for t\in[0,T]\right\rbrace, 
\end{equation*}
which is an unbounded open set on $\RR\times X$. Observe that $F:\RR\times X\to X$ is continuous and maps bounded subsets of $\RR\times X$ whose closure is contained on $\Lambda$ into relative compact sets of $X$ (i.e., $F:\RR\times X\to X$ is completely continuous on $\Lambda$). So that the Leray-Schauder degree theory applies for this kind of bounded open subsets of $\Lambda$.

To state the main result of this section it will be convenient to introduce some notation. Given $\Delta\subseteq\Lambda$ and $\beta\in\RR$, we denote by $\Delta_{\beta}$ the vertical section $\Delta_{\beta}=\lbrace \xi\in X:(\beta,\xi)\in \Delta\rbrace$. From now on, without loss of generality the following identification
\begin{equation*}
(\beta,x,y,f)\equiv \rho(t):=a+x+\beta\int_0^t f(s)ds\for t\in [0,T]
\end{equation*}
can be used to simplify the notation, taking into account that if $(\beta,x,y,f)\in\Lambda$ is a solution to \eqref{2.3}-\eqref{2.4} then $\rho(t)$ is a positive $T-$periodic solution of \eqref{2.2}.

We are now in the position to prove, by using continuation arguments of degree theory, that \eqref{2.3}-\eqref{2.4} has many solutions.

\begin{Lem}\label{lema 2}
There exist a positive constant $\beta_*>0$ and a connected set $\mathcal{C}$ of solutions to \eqref{2.3}-\eqref{2.4} satisfying
\begin{equation}\label{nn2}
\lbrace\beta: (\beta,x,y,f)\in\mathcal{C}\rbrace\supseteq (0,\beta_*]
\end{equation}
and one of the following two possibilities:
\begin{description}
\item[(i)] $\mathcal{C}$ is unbounded;
\item[(ii)] $\inf\lbrace m_{\rho}:(\beta,x,y,f)\in\mathcal{C}\rbrace$=0.
\end{description}
\end{Lem}
\begin{proof}
According to Notation \ref{nota1} we choose $0<\varepsilon<a$ such that
\begin{equation}\label{2n}
\varepsilon<\sqrt{\lambda\alpha}<a-\varepsilon.
\end{equation}
Consider the set
\begin{equation*}
\Omega:=\left\lbrace(\beta,x,y,f):|\beta|<\frac{\varepsilon}{2T(y_0+ \delta_0\|\sigma\|_{\infty})},|x|<a-\varepsilon,|y|<y_0,\|f\|_{\infty}<y_0+\delta_0\|\sigma\|_{\infty}\right\rbrace, 
\end{equation*}
where
\begin{equation*}
y_0=\frac{(1+\lambda)|\overline{\sigma}|}{2\sqrt{\lambda\alpha}}+\varepsilon,\qquad \delta_0:=\frac{1+\lambda}{2\sqrt{\lambda\alpha}};
\end{equation*}
and observe that $\Omega$ is a bounded and open set on $\RR\times X$ whose closure is contained in $\Lambda$. Here, we shall prove that $F[0,z]\neq z$ for all $z\in\partial\Omega_0$. Indeed, this set can be divided into three (nondisjoint) subsets:
\begin{equation}\label{nn1}
\partial\Omega_0=A_1\cup A_2\cup A_3,
\end{equation}
where $A_1$ is a subset of elements $(x,y,f)\in X$ such that $|x|=a-\varepsilon$, $A_2$ is composed by elements $(x,y,f)\in X$ for which $|y|=y_0$ and $A_3$ is a subset of $X$ whose elements verify that $\|f\|_{\infty}=y_0+\delta_0\|\sigma\|_{\infty}$. On the other hand if we assume that $F[0,\cdot]$ has a fixed point $(x,y,f)\in\partial\Omega_0$, then easily one checks that it should lie on the finite dimensional subspace $\RR^2\times\mbox{span}\lbrace 1,\sigma\rbrace$ (notice that $F[0,\Lambda]\subseteq \RR^2\times\mbox{span}\lbrace 1,\sigma\rbrace$).  In particular this implies that
\begin{equation*}
f(t)=\gamma+\delta\sigma(t),\for t\in [0,T],
\end{equation*}
for some $\gamma,\delta\in\RR$ verifying \eqref{n1}. It means that
\begin{equation}\label{*}
\left.
\begin{array}{rcl}
\ds\int_0^T(\gamma+\delta\sigma(s))ds&=&0\\\\
\ds\frac{1}{T}\int_0^T(\gamma+\delta\sigma(s))^2ds&=&\frac{(1+\lambda)^2}{4\lambda}\\
\gamma&=&y\\
\delta&=&\frac{1+\lambda}{2(a+x)}
\end{array}
\right\}.
\end{equation}
Solving the system we obtain only one solution on $\Lambda_0$:
\begin{equation}\label{n3}
\delta^*=\delta_0,\qquad y^*=\gamma^*=-\overline{\sigma}\delta_0,\qquad x^*=\sqrt{\lambda\alpha}-a.
\end{equation}
Now, the fact that $F[0,\cdot]$ does not have fixed points on $\partial\Omega_0$ is a direct consequence of \eqref{2n}, \eqref{nn1} and the definitions of $y_0$ and $\delta_0$.

Having proved that the Leray-Schauder degree of $(I-F)[0,\cdot]$ is well defined on the section $\Omega_0$, we need to compute it. Since $F[\cdot,\overline{\Omega}_0]$ is included in a $4-$dimensional subspace, according to \cite[Theorem 8.7]{DEIM}, we obtain that
\begin{equation}\label{n4}
d_{LS}(I-F[0,\cdot],\Omega_0,0)=d_{B}(I_{\RR^4}-F[0,\cdot]_{\RR^4},\Omega_0\cap\RR^4,0),
\end{equation}
where $\RR^2\times\mbox{span}\lbrace 1,\sigma\rbrace$ is identified by $\RR^4$, 
\begin{multline*}
(I_{\RR^4}-F[0,\cdot]_{\RR^4})(x,y,\gamma,\delta)=\left(\int_0^T(\gamma+\delta\sigma(s))ds,\right.\\ \left. \frac{1}{T}\int_0^T(\gamma+\delta\sigma(s))^2ds-\frac{(1+\lambda)^2}{4\lambda},\gamma-y,\delta-\frac{1+\lambda}{2(a+x)} \right)
\end{multline*}
and $\Omega_0\cap \RR^4$ is an open set containing the unique solution of \eqref{*} given by \eqref{n3}. By using Brouwer degree theory we have that
\begin{equation}\label{n5}
\mbox{sgn } \mathbb{J}=d_{B}(I_{\RR^4}-F[0,\cdot]_{\RR^4},\Omega_0\cap\RR^4,0),
\end{equation}
where $\mathbb{J}$ denotes the determinant of the Jacobian matrix associated to the function $I_{\RR^4}-F[0,\cdot]_{\RR^4}$ at the point $(x^*,y^*,\delta^*,\gamma^*)$ (given by \eqref{n3}); i.e.,
\begin{equation*}
\begin{array}{rcl}
\mathbb{J}&=&\begin{vmatrix}
0  &  0   & T                &                              \ds\int_0^T\sigma(s)ds\\
0  &  0   & \ds\frac{2}{T}\int_0^T(\gamma^*+\delta^*\sigma(s))ds & \ds\frac{2}{T}\int_0^T(\gamma^*+\delta^*\sigma(s))\sigma(s)ds\\
0  &  -1  & 1       &  0\\
\ds\frac{1+\lambda}{2(a+x^*)^2}  & 0  & 0 & 1
\end{vmatrix}\\\\
  &=&\ds\frac{1+\lambda}{(a+x^*)^2}\int_0^T(\gamma^*+\delta^*\sigma(s))\sigma(s)ds= \ds\frac{(1+\lambda)\delta^*}{(a+x^*)^2}\int_0^T(\sigma(s)-\overline{\sigma})\sigma(s)ds.
\end{array}
\end{equation*}
Observe that, since $\alpha>0$ then $\mathbb{J}\neq 0$. Combining \eqref{n4} with \eqref{n5} we see that, indeed, $d_{LS}(I-F[0,\cdot],\Omega_0,0)\neq 0$.

Under these conditions the classical Leray-Schauder continuation theorem (see, e.g. \cite{LerSch} and also \cite{Mtmna,MRZccm}) provides the existence of a connected set $\mathcal{C}$ composed by solutions to \eqref{2.2}-\eqref{2.3} verifying \eqref{nn2} for some $\beta_*>0$ and satisfying one of the following three possibilities (see Figure 1): (a) $\mathcal{C}$ is unbounded on $\Lambda$, (b) $\mbox{dist}(\mathcal{C},\partial\Lambda)=0$ or (c) $\mathcal{C}_0\cap (\Lambda_0\setminus\overline{\Omega}_0)\neq\emptyset$. Concerning to the possibility (c), one observes that \eqref{*} has only one solution on $\Lambda_0$, which must belong to $\Omega_0$, obtaining a contradiction. On the other hand, the possibilities (a) and (b) imply (i) and (ii), completing the proof.
\begin{figure}[!h]
\begin{center}
\includegraphics[scale=1.2]{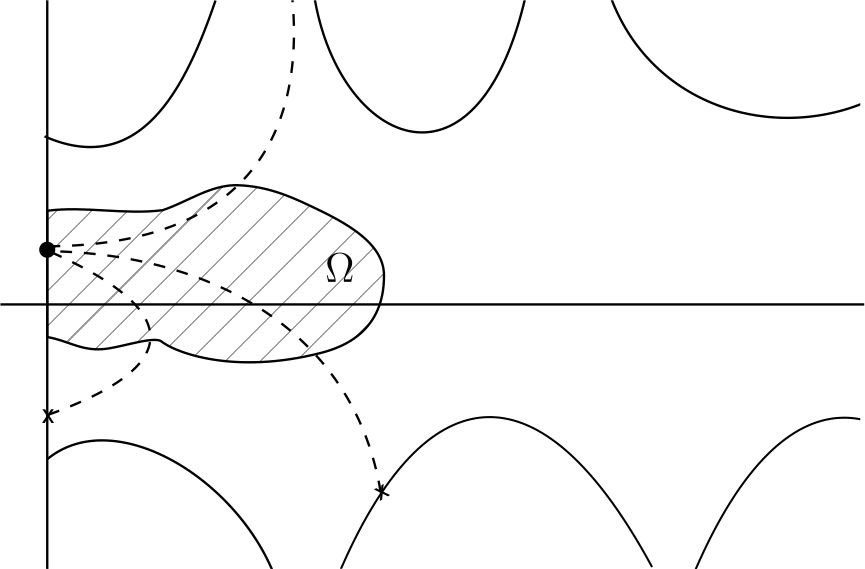}\put(-246,66){\small $0$}\put(-13,66){\small $\beta$}\put(-164, 124){\small (a)}\put(-247, 155){\small $X$}\put(-143,51){\small (b)\put(-78,-5){\small (c)}}\end{center}\caption{The possibilities of the set $\mathcal{C}$}
\end{figure}
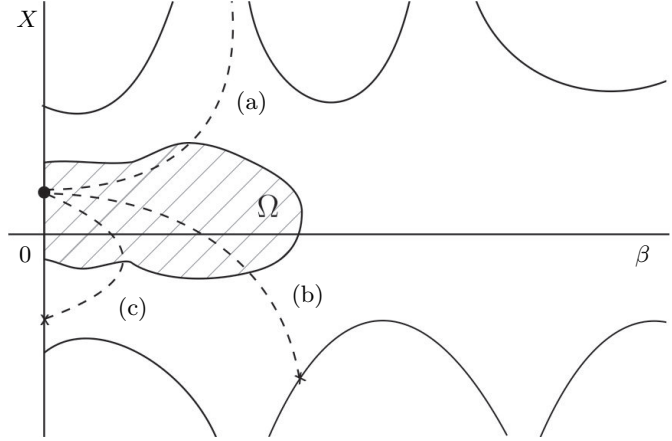
\end{proof}

\section{A priori bounds on the components of $\mathcal{C}$}\label{sec4}

In this section we will analyse in detail the behaviour of the set $\mathcal{C}$ in order to determinate when either it intersects $\partial\Lambda$ (i.e. (ii) holds) or it connects with points coming from infinity (i.e. $\beta_*=+\infty$). With this aim we are going to assume that there exists $\rho_0>0$ such that
\begin{equation}\label{3.1}
\rho(t)\geq\rho_0\for t\in [0,T]
\end{equation}
for any $(\beta,x,y,f)\in\mathcal{C}$. Under this condition we shall prove the boundedness of $\mathcal{C}$ over its components $x, y$ and $f$.
\begin{Lem}\label{lema 3.1}
The inequality
\begin{equation*}
\|f\|_{\infty}\leq \frac{1+\lambda}{2\rho_0}\left(\frac{\beta(1+3\lambda)T}{2\lambda}+\|\widetilde{h}\|_1\right) 
\end{equation*}
holds for any $(\beta,x,y,f)\in \mathcal{C}$.
\end{Lem}
\begin{proof}
We let $(\beta,x,y,f)\in\mathcal{C}$ such that \eqref{3.1} holds. From the second equation of \eqref{2.3} we notice that $f$ is absolutely continuous on the interval $[0,T]$ and 
\begin{equation*}
f'(t)\left[a+x+\beta\int_0^t f(s)ds\right]=\left(\frac{\lambda-1}{\lambda+1}\right)\beta f^2(t)+\left(\frac{1+\lambda}{2}\right)\left(-\beta+\widetilde{h}(t)\right).   
\end{equation*}
According to the first equation of \eqref{2.4} and by using \eqref{3.1} we prove that
\begin{equation*}
\rho_0\|f'\|_1\leq\frac{1+\lambda}{2}\left(\frac{\beta(1+3\lambda)T}{2\lambda}+\|\widetilde{h}\|_1\right). 
\end{equation*}
We conclude the proof taking into account that $f$ changes its sign (note that $\overline{f}=0$), therefore $\|f\|_{\infty}\leq\|f'\|_1$.
\end{proof}
\begin{Rem}\label{remark 3.1}
We note that
\begin{equation*}
|y|=|f(0)|\leq \frac{1+\lambda}{2\rho_0}\left(\frac{\beta(1+3\lambda)T}{2\lambda}+\|\widetilde{h}\|_1\right)
\end{equation*}
for any $(\beta,x,y,f)\in\mathcal{C}$.
\end{Rem}
\begin{Lem}\label{lema 3.2}
The inequality
\begin{equation*}
|x|\leq a+\frac{\beta T(1+\lambda)}{2\sqrt{\lambda}}\left(1+\frac{\widetilde{H}^{\frac{1+\lambda}{2\lambda}}}{(\beta T+\widetilde{H})^{\frac{1+\lambda}{2\lambda}}-\widetilde{H}^{\frac{1+\lambda}{2\lambda}}}\right) 
\end{equation*}
holds for any $(\beta,x,y,f)\in\mathcal{C}$.
\end{Lem}
\begin{proof}
We let $(\beta,x,y,f)\in\mathcal{C}$ such that \eqref{3.1} holds.  By using the Cauchy--Bunyakovsky--Schwarz inequality we have
\begin{equation}\label{3.1'}
x\geq -\left(a+\frac{\beta(1+\lambda)T}{2\sqrt{\lambda}}\right). 
\end{equation}
Since
\begin{equation*}
M_{\rho}-m_{\rho}\leq\beta\int_0^T|f(s)|ds\leq\frac{\beta(1+\lambda)T}{2\sqrt{\lambda}},
\end{equation*}
then
\begin{equation}\label{3.2}
x\leq m_{\rho}+\frac{\beta(1+\lambda)T}{2\sqrt{\lambda}} 
\end{equation}
(notice that $\rho(0)=a+x$). 

Consider the function $u(t)=\rho^{2/(1+\lambda)}(t)$, by \eqref{2.1.1} we know that $u$ is a $T-$periodic solution to \eqref{2.1} such that
\begin{equation*}
m_{u}=m_{\rho}^{\frac{2}{1+\lambda}},\qquad
M_u=M_{\rho}^{\frac{2}{1+\lambda}}\leq \left[m_{\rho}+\frac{\beta(1+\lambda)T}{2\sqrt{\lambda}}\right]^{\frac{2}{1+\lambda}}. 
\end{equation*}
Taking into account the latter inequalities, we integrate both sides of \eqref{2.1} over $[0,T]$ obtaining
\begin{equation*}
\frac{\beta T+\widetilde{H}}{\left(m_{\rho}+\frac{\beta(1+\lambda)T}{2\sqrt{\lambda}}\right)^{\frac{2\lambda}{1+\lambda}}}\leq\frac{\widetilde{H}}{m_{\rho}^{\frac{2\lambda}{1+\lambda}}}.
\end{equation*}
So that
\begin{equation}\label{3.3}
m_{\rho}\leq\frac{\beta T(1+\lambda)\widetilde{H}^{\frac{1+\lambda}{2\lambda}}}{2\sqrt{\lambda}\left[(\beta T+\widetilde{H})^{\frac{1+\lambda}{2\lambda}}-\widetilde{H}^{\frac{1+\lambda}{2\lambda}}\right]}.
\end{equation}
Now, the result immediately follows from \eqref{3.1'}, \eqref{3.2} and \eqref{3.3}.
\end{proof}

\begin{Pro}\label{prop 1}
Assume that (ii) does not hold. Then
\begin{equation*}
\mathcal{C}_{\beta}\neq\emptyset\for\beta>0.
\end{equation*}
\end{Pro}
\begin{proof}
Using a contradiction argument, assume that
\begin{equation*}
\left\lbrace\beta:(\beta,x,y,f)\in\mathcal{C}\right\rbrace\subseteq (0,\beta^*] 
\end{equation*}
for some $\beta^*>0$. According to Lemma \ref{lema 3.1} and Remark \ref{remark 3.1} we prove the boundedness of the components $f$ and $y$ of the connected set $\mathcal{C}$. Since
\begin{equation*}
\lim_{\beta\to 0^+}\frac{\beta (1+\lambda)T}{2\sqrt{\lambda}}\left(1+\frac{\widetilde{H}^{\frac{1+\lambda}{2\lambda}}}{(\beta T+\widetilde{H})^{\frac{1+\lambda}{2\lambda}}-\widetilde{H}^{\frac{1+\lambda}{2\lambda}}}\right)=\sqrt{\lambda}\widetilde{H}, 
\end{equation*} 
then by Lemma \ref{lema 3.2},
\begin{equation*}
|x|\leq\max_{\beta\in [0,\beta^*]}\left[ a+\frac{\beta (1+\lambda)T}{2\sqrt{\lambda}}\left(1+\frac{\widetilde{H}^{\frac{1+\lambda}{2\lambda}}}{(\beta T+\widetilde{H})^{\frac{1+\lambda}{2\lambda}}-\widetilde{H}^{\frac{1+\lambda}{2\lambda}}}\right)\right]. 
\end{equation*}
Therefore $\mathcal{C}$ is bounded, obtaining a contradiction with Lemma \ref{lema 2}. 
\end{proof}

\section{Proof of the main results}\label{sec5}
The proof of Theorem \ref{teo1} is especially simple after combining Lemma \ref{lema 2}, Proposition \ref{prop 2.1} and the equivalence between $T-$periodic solutions of \eqref{2.1} and positive $T-$periodic solutions of \eqref{2.2} discussed in Section 2.
\begin{proof}[Proof of Theorem \ref{teo1}]
From Lemma \ref{lema 2}, we define $\beta_*>0$ and consider $\mathcal{C}$ the connected set of solutions to \eqref{2.3}-\eqref{2.4} satisfying \eqref{nn2}. Concerning to this solutions, for any $\beta\in (0,\beta_*]$ we have $\rho(t)=(\beta,x,y,f)$ a positive $T-$periodic solution of \eqref{2.2}. In view of \eqref{2.1.1}, $u$ is also a $T-$periodic solution of \eqref{2.1}. The result is obtained observing that $v(t)=\beta^{-\frac{1}{1+\lambda}}u(t)$ for all $t\in [0,T]$
is a $T-$periodic solution of \eqref{2}.
\end{proof}

The remaining of this section is devoted to prove Corollary \ref{cor 1}. To this end we will consider $\widetilde{h}(t):=h(t)-\overline{h}\in\widetilde{L}(\RR/T\ZZ)$
and we will apply Theorem \ref{teo1} to obtain the existence of a $T-$periodic solution of \eqref{2} provided that $\beta\in (0,\beta_*]$. This perturbation result may give global information on the existence of $T-$periodic solutions to \eqref{1} if a positive lower bound of $\beta_*>0$ is obtained. Taking into account that $\beta_*>0$ is related to a connected set $\mathcal{C}$ of solutions to \eqref{2.2}, our strategy shall consist in analizing the properties of $\mathcal{C}$ to determine an estimate for the parameter $\beta_*>0$. For instance, in view of Proposition \ref{prop 1} we have two possibilities: either $\beta_*=+\infty$ or 
\begin{equation}\label{5.1}
\inf\left\lbrace m_{\rho}:(\beta,x,y,f)\in\mathcal{C}\right\rbrace=0. 
\end{equation}
If $\beta_*=+\infty$, it is simple to conclude that \eqref{1} has a $T-$periodic solution. Hence, there is no loss of generality in assuming \eqref{5.1} in order to estimate $\beta_*>0$. Unfortunately, the general case in which $\widetilde{h}\in\widetilde{L}(\RR/T\ZZ)$ is a nontrivial function presents some difficulties which seem to require a deeper treatment. Nevertheless the situation is more simple if the forcing term $\widetilde{h}$ under consideration is piecewise-constant. 

We consider the simplest case as a starting point in order to avoid getting lost in secondary details and keeping the exposition at a rather simple level. However, the arguments of this section can be extended to more general situations by making small changes (see Section \ref{sec6}). Therefore this section is devoted to consider 
\begin{equation*}
\widetilde{h}(t)=
\begin{cases}
h_*,&\quad {\mbox { if } } t\in [0,T/2),\\
-h_*,&\quad {\mbox { if } }t\in [T/2,T),
\end{cases}
\end{equation*}
where $h_*>0$. Under this framework we establish the following result.
\begin{Pro}\label{pro 3}
The following estimation
\begin{equation*}
\beta_*\geq\frac{\lambda(1-\lambda)h_*}{4(1+\lambda)^2}
\end{equation*}
holds.
\end{Pro}
\begin{proof}
Observe that \eqref{2.1} can be studied as two alternating autonomous equations:
\begin{equation}\label{5.2}
\begin{cases}
u''=\frac{-\beta^2+\beta h_*}{u^{\lambda}},&\quad {\mbox { if } } t\in [0,T/2),\\\\
u''=\frac{-\beta^2-\beta h_*}{u^{\lambda}},&\quad {\mbox { if } }t\in [T/2,T).
\end{cases}
\end{equation}
Without loss of generality we can assume \eqref{5.1}. Using a contradiction argument,  we assume that
\begin{equation*}
\beta_*<\frac{\lambda(1-\lambda)h_*}{4(1+\lambda)^2}.
\end{equation*}
In view of \eqref{5.1}, we can choose $\rho(t)=(\beta,x,y,f)\in\mathcal{C}$ with $\beta\in (0,\beta_*]$ a positive $T-$periodic solution of \eqref{2.2} such that
\begin{equation}\label{5.3}
\beta_*<\frac{\lambda(1-\lambda)h_*}{4(1+\lambda)^2}-\frac{4m_{\rho}\sqrt{\lambda}}{T(1+\lambda)}.
\end{equation}
Consider $u$ the $T-$periodic solution of \eqref{2.1} given by \eqref{2.1.1}. We already know that 
\begin{gather}
\notag 
m_u=m_{\rho}^{\frac{2}{1+\lambda}},\\
\label{5.4'} 
M_u=M_{\rho}^{\frac{2}{1+\lambda}}\leq \left[m_{\rho}+\frac{\beta(1+\lambda)T}{2\sqrt{\lambda}}\right]^{\frac{2}{1+\lambda}}. 
\end{gather}
Moreover, according to Lemmas \ref{lemma 4.1} and \ref{lemma 4.2} (see the appendix, i.e., Section \ref{sec7}), multiplying both sides of the first equation of \eqref{5.2} by $u'$, the integration over the interval $[T/4,t]$ combined with the fact that $u'(t)>0$ for $[T/4,T/2]$ yields
\begin{equation}\label{5.4''}
\int_{m_{\rho}^{\frac{2}{1+\lambda}}}^{u(T/2)}\left[s^{1-\lambda}-m_{\rho}^{\frac{2(1-\lambda)}{1+\lambda}}\right]^{-1/2}ds=\frac{T}{4}\sqrt{\frac{2(-\beta^2+\beta h_*)}{1-\lambda}}.
\end{equation}
Since $u''> 0$ on $[0,T/2)$, then $\max_{t\in [0,T/2]}u(t)=u(T/2)$ and
\begin{multline*}
\int_{m_{\rho}^{\frac{2}{1+\lambda}}}^{u(T/2)}\left[ s^{1-\lambda}-m_{\rho}^{\frac{2(1-\lambda)}{1+\lambda}}\right]^{-1/2}ds\\=\int_{m_{\rho}^{\frac{2}{1+\lambda}}}^{u(T/2)}\frac{s^{\lambda}\left[s^{1-\lambda}-m_{\rho}^{\frac{2(1-\lambda)}{1+\lambda}}\right]^{-1/2}(1-\lambda)s^{-\lambda}}{1-\lambda}ds\\
\leq\frac{u^{\lambda}(T/2)}{1-\lambda}\int_{m_{\rho}^{\frac{2}{1+\lambda}}}^{u(T/2)}\left[s^{1-\lambda}-m_{\rho}^{\frac{2(1-\lambda)}{1+\lambda}}\right]^{-1/2}(1-\lambda)s^{-\lambda}ds\\
=\frac{2u^{\lambda}(T/2)}{1-\lambda}\left[u^{1-\lambda}\left(T/2\right)-m_{\rho}^{\frac{2(1-\lambda)}{1+\lambda}}\right]^{1/2}.
\end{multline*}
Consequently,
\begin{equation*}
\frac{u^{\lambda}(T/2)}{1-\lambda}\left[u^{1-\lambda}(T/2)-m_{\rho}^{\frac{2(1-\lambda)}{1+\lambda}}\right]^{1/2}\geq\frac{T}{8}\sqrt{\frac{2(-\beta^2+\beta h_*)}{1-\lambda}}. 
\end{equation*}
Squaring both sides of the above inequality and taking in mind that $u(T/2)\geq m_{\rho}^{2/(1+\lambda)}$, one has
\begin{equation}\label{5.4}
u^{1+\lambda}\left(T/2\right)\geq m_{\rho}^2+\frac{(1-\lambda)(-\beta^2+\beta h_*)T^2}{32}.
\end{equation}  

Arguing analogously as before on the interval $[T/2,3T/4]$, with respect to the second equation of \eqref{5.2} we arrive to
\begin{equation*}
M_u^{1+\lambda}-M_u^{2\lambda}u^{1-\lambda}(T/2)\geq\frac{T^2(\beta^2+\beta h_*)(1-\lambda)}{32}.
\end{equation*}
From \eqref{5.4} it follows that
\begin{equation*}
M_u^{1+\lambda}\geq m_{\rho}^2+\frac{T^2(1-\lambda)\beta h_*}{16}.
\end{equation*}
In view of \eqref{5.4'},
\begin{equation*}
\left[m_{\rho}+\frac{\beta(1+\lambda)T}{2\sqrt{\lambda}}\right]^2\geq m_{\rho}^2+\frac{T^2(1-\lambda)\beta h_*}{16}, 
\end{equation*}
and, consequently,
\begin{equation*}
\frac{\beta (1+\lambda)^2 T}{4\lambda}+\frac{(1+\lambda) m_{\rho}}{\sqrt{\lambda}}\geq\frac{T(1-\lambda)h_*}{16}.
\end{equation*}  
Solving $\beta$ in the last inequality one obtains a contradiction with \eqref{5.3}. The proof is thus complete.
\end{proof}

We are now in the position to prove Corollary \ref{cor 1}.
\begin{proof}[Proof of Corollary \ref{cor 1}]
Define
\begin{equation*}
\widetilde{h}(t)=
\begin{cases}
\frac{h_1+h_2}{2},&\quad {\mbox { if } } t\in [0,T/2),\\
-\frac{h_1+h_2}{2},&\quad {\mbox { if } }t\in [T/2,T).
\end{cases}
\end{equation*}
Combining Theorem \ref{teo1} and Proposition \ref{pro 3}, the condition \eqref{hip} implies that \eqref{2} has a $T-$periodic solution provided $\beta:=-\overline{h}=(h_2-h_1)/2$. Since $h(t)=\overline{h}+\widetilde{h}(t)$ for all $t\in [0,T]$, the proof is complete.
\end{proof}

\section{The standard extension to general piecewise-constant forcing terms}\label{sec6}

We already pointed out that it is easy to extend Corollary \ref{cor 1} by considering  more general forcing terms and doing slight modifications to the proof of Proposition \ref{pro 3}. The details of this section can be easily verified; therefore, the proofs  will be omitted.

An important aspect to note now is that Lemmas \ref{lemma 4.1} and \ref{lemma 4.2} can be used to extend directly Corollary \ref{cor 1} to forcing terms having two asymmetrical weights. Hence, if we consider
\begin{equation*}
h(t)=
\begin{cases}
h_1,&\quad {\mbox { if } } t\in [0,\eta),\\
-h_2,&\quad {\mbox { if } }t\in [\eta,T),
\end{cases}
\end{equation*}
where $\eta\in (0,T)$ and $h_1, h_2$ are positive constants; it should not be too difficult for the reader to prove 
\begin{Cor}\label{cor 3}
Assuming that
\begin{equation*}
0<\frac{h_2(T-\eta)-h_1\eta}{h_2+h_1}\leq\frac{\lambda(1-\lambda)\eta(T-\eta)}{2T(1+\lambda)^2-\lambda(1-\lambda)(T-2\eta)},
\end{equation*}
then \eqref{1} has a $T-$periodic solution.
\end{Cor}
To prove this result one has to estimate the parameter $\beta_*>0$ obtained by Theorem \ref{teo1} when the forcing term under consideration is
\begin{equation*}
\widetilde{h}(t)=
\begin{cases}
h_1^*,&\quad {\mbox { if } } t\in [0,\eta),\\
-h_2^*,&\quad {\mbox { if } }t\in [\eta,T),
\end{cases}
\end{equation*}
where now $h_1^*, h_2^*$ are positive constants verifying that $h_1^*\eta=h_2^*(T-\eta)$. Arguing as in Proposition \ref{pro 3}, under this framework it is easy to see that
\begin{equation*}
\beta_*\geq\frac{\lambda(1-\lambda)h_1^*\eta}{2T(1+\lambda)^2-\lambda(1-\lambda)(T-2\eta)}.
\end{equation*}

Now, it is the right time to make a more in-depth analysis of the proof of Proposition \ref{pro 3} and understanding the applicability of our result to more general class of piecewise-constant functions. Observe that until now we have worked with piecewise-constant forcing terms with only two values in order to exploit the fact that the $T-$periodic solutions of these equations attain their extremal values at particular points ($t_m=\eta/2, t_M=(\eta+T)/2$). However this property only seems crucial for optimizing the estimation of the parameter $\beta_*>0$ provided in Theorem \ref{teo1}, but it does not play a fundamental role if the aim is to obtain some estimation of $\beta_*>0$, adapting the proof of Proposition \ref{pro 3}. For instance, if we consider 
\begin{equation*}
h(t)=h_i,\quad{\mbox { if } } t\in \left[\frac{(i-1)T}{n},\right.\left.\frac{iT}{n}\right),
\end{equation*}
where $h_i\in\RR$ for each $i=1,\ldots, n$ ($n$ is the number of weights). In this setting one can verify
\begin{Cor}\label{cor 2}
Assuming that
\begin{equation}\label{hip 2}
0<\frac{-\sum_{i=1}^{n}h_i}{\min\left\lbrace  nh_j:nh_j>\sum_{i=1}^{n}h_i\right\rbrace-\sum_{i=1}^{n}h_i}\leq\frac{\lambda(1-\lambda)}{2n^2(1+\lambda)^2+\lambda(1-\lambda)},
\end{equation}
then \eqref{1} has a $T-$periodic solution.
\end{Cor}
Observe that the second inequality of \eqref{hip 2} only can be satisfied if $h_i>0$ for some $i=1,\ldots,n$; which is at same time a necessary condition to guarantee the existence of $T-$periodic solutions of \eqref{1}. Again, the key to prove this result consists in estimating the parameter $\beta_*>0$ obtained in Theorem \ref{teo1} be considering now the function
\begin{equation*}
\widetilde{h}(t)=h_i^*,\quad{\mbox { if } } t\in \left[\frac{(i-1)T}{n},\right.\left.\frac{iT}{n}\right),
\end{equation*}
where $h_i^*\in\RR$ for each $i=1,\ldots, n$ and $\sum_{i=1}^n h_i^*=0$. Under this framework and after doing slight modifications to the proof of Proposition \ref{pro 3}, we can prove that
\begin{equation}\label{6.2}
\beta_*\geq\frac{\lambda(1-\lambda)\min\lbrace h_i^*:h_i^*>0\rbrace}{2n^2(1+\lambda)^2+\lambda(1-\lambda)}.
\end{equation}

Indeed, we let $t_m\in [0,T]$ a point where the function $u$ (given by Proposition \ref{pro 3}) attains its minimum value. Without loss of generality we can assume that $t_m\in I_i$, where $I_i:=[a_i,b_i]$ is an interval such that $\widetilde{h}|_{I_i}\equiv h_i^*>0$. Following the argument of Proposition \ref{pro 3} we can obtain
\begin{equation*}
\int_{m_{\rho}^{\frac{2}{1+\lambda}}}^{u(\vartheta)}\left[s^{1-\lambda}-m_{\rho}^{\frac{2(1-\lambda)}{1+\lambda}}\right]^{-1/2}ds=\left.|\vartheta-t_m\right.|\sqrt{\frac{2(-\beta^2+\beta h_i^*)}{1-\lambda}}
\end{equation*}
instead of \eqref{5.4''}, where $\vartheta\in \lbrace a_i,b_i\rbrace$ is defined such that
\begin{equation*}
|\vartheta-t_m|=\max\lbrace t_m-a_i,b_i-t_m\rbrace\geq\frac{T}{2n}.
\end{equation*}
Consequently, 
\begin{equation*}
\int_{m_{\rho}^{\frac{2}{1+\lambda}}}^{u(\vartheta)}\left[s^{1-\lambda}-m_{\rho}^{\frac{2(1-\lambda)}{1+\lambda}}\right]^{-1/2}ds\geq\frac{T}{2n}\sqrt{\frac{2(-\beta^2+\beta h_i^*)}{1-\lambda}}.
\end{equation*}
Arguing analogously as in Proposition \ref{pro 3} we arrive to
\begin{equation*}
u^{1+\lambda}(\vartheta)-m_{\rho}^{\frac{2(1-\lambda)}{1+\lambda}}u^{2\lambda}(\vartheta)\geq\frac{T^2(1-\lambda)(-\beta^2+\beta h_i^*)}{8n^2}.
\end{equation*}  
Since $u(\vartheta)\geq m_{\rho}^{2/(\lambda+1)}$, \eqref{5.4'} yields
\begin{equation*}
\left[m_{\rho}+\frac{\beta(1+\lambda)T}{2\sqrt{\lambda}}\right]^{2}\geq m_{\rho}^2+\frac{T^2(1-\lambda)(-\beta^2+\beta h_i^*)}{8n^2}. 
\end{equation*}
Solving $\beta$ and taking into account that $m_{\rho}$ can became arbitrarily small (see \eqref{5.1}) one proves the desired estimation. 
\begin{Rem}\label{rema1}
Of course, after introducing more notation one can easily prove similar results dealing with asymmetrical piecewise-constant forcing terms having a finite number of pieces. This framework will be omitted here. 
\end{Rem} 
\section{Apendix: Two auxiliary lemmas}\label{sec7}
This section is dedicated to the study of equation \eqref{1} when the forcing term under consideration is piecewise-constant with only two values, i.e.,
\begin{equation*}
h(t)=
\begin{cases}
h_1,&\quad {\mbox { if } } t\in [0,\eta),\\
-h_2,&\quad {\mbox { if } }t\in [\eta,T),
\end{cases}
\end{equation*}
where $\eta\in (0,T)$ and $h_1, h_2$ are positive constants. The goal of this section is to obtain an explicit description of the $T-$periodic solutions of \eqref{1} and take advantage of this knowledge to optimize some results of this paper dealing with this class of functions. We point out that the results of this appendix are not completely new nor play a crucial role in the proofs of the main results of this note (see Section \ref{sec6}). That is why we dedicate this special section in this paper. The proofs of Lemmas \ref{lemma 4.1} and \ref{lemma 4.2} can be found in \cite{HZpreprint} under a different framework, but the same argument works.

From now on, for simplicity, we shall assume the previous setting with $\lambda\in (0,1)$. The following lemma establishes that if $u$ is a $T-$periodic solution of \eqref{1} then its value at $\eta$ is necessary equal to its value at the extremes points on the interval $[0,T]$.
\begin{Lem}\label{lemma 4.1}
If $u$ is a $T-$periodic solution of \eqref{1} then
\begin{equation*}
u(0)=u\left(\eta\right)=u(T).
\end{equation*}
\end{Lem}
\begin{proof}
According to sign properties of the function $h$ there exist $t_m\in (0,\eta)$ and $t_M\in (\eta,T)$ such that $m_u=u(t_m)$ and $M_u=u(t_M)$. Multiplying both sides of \eqref{1} by $u'$ and integrating on $[t_m,\eta]$ we obtain
\begin{equation}\label{7.1}
u'^2(\eta)=\frac{2h_1\left(u^{1-\lambda}(\eta)-m_u^{1-\lambda}\right)}{1-\lambda}. 
\end{equation}
Further, multiplying both sides of \eqref{1} by $u'$ and integrating on $[\eta, t_M]$ we get
\begin{equation}\label{7.2}
u'^2(\eta)=\frac{2h_2\left(-u^{1-\lambda}(\eta)+M_u^{1-\lambda}\right)}{1-\lambda}.
\end{equation}
In view of \eqref{7.1} and \eqref{7.2} we conclude that
\begin{equation*}
u^{1-\lambda}(\eta)=\frac{h_1 m_u^{1-\lambda}+h_2 M_u^{1-\lambda}}{h_1+h_2}.
\end{equation*}
Arguing analogously as before on the intervals $[0,t_m]$ and $[t_M,T]$ respectively, and using the periodicity conditions we can prove that
\begin{equation}\label{7.3}
h_1\left(u^{1-\lambda}(0)-m_u^{1-\lambda}\right)=h_2\left(M_u^{1-\lambda}-u^{1-\lambda}(T)\right).  
\end{equation} 
From \eqref{7.3} it follows that $u^{1-\lambda}(0)=u^{1-\lambda}(\eta)$, concluding the proof.
\end{proof}
As a consequence of the previous lemma we get that each $T-$periodic solution of \eqref{1} attains its maximum and minimum value at the points $\eta/2$ and $(\eta+T)/2$, respectively.
\begin{Lem}\label{lemma 4.2}
If $u$ is a $T-$periodic solution of \eqref{1} then
\begin{equation*}
m_u=u\left(\frac{\eta}{2}\right) ,\qquad M_u=u\left( \frac{\eta+T}{2}\right).
\end{equation*}
\end{Lem}
\begin{proof}
We shall prove that $t_m=\eta/2$. Indeed, since $u''>0$ for $t\in [0,t_m)$ then $u'(t)<0$ for $t\in [0,t_m)$ (note that $u'(t_m)=0$). Multiplying both sides of \eqref{1} by $u'$ and integrating over $[t,t_m]$, we arrive to
\begin{equation*}
u'(t)=-\sqrt{\frac{2h_1}{1-\lambda}}\sqrt{u^{1-\lambda}(t)-m_u^{1-\lambda}}\for t\in [0,t_m],
\end{equation*}
and, consequently,
\begin{equation*}
\int_{m_u}^{u(0)}\frac{ds}{\sqrt{s^{1-\lambda}-m_u^{1-\lambda}}}=t_m\sqrt{\frac{2h_1}{1-\lambda}}.
\end{equation*}
Analogously, the double integration over the interval $[t_m,t]$, with respect to the fact that $u'(t)>0$ for $t\in (t_m,\eta]$, yields
\begin{equation*}
\int_{m_u}^{u(\eta)}\frac{ds}{\sqrt{s^{1-\lambda}-m_u^{1-\lambda}}}=(\eta-t_m)\sqrt{\frac{2h_1}{1-\lambda}}.
\end{equation*}
According to Lemma \ref{lemma 4.1}, by subtracting the above-mentioned identities we obtain that $t_m=\eta/2$.

The proof of the relation $t_M=(\eta+T)/2$ is analogous and it will be omitted. 
\end{proof}
To conclude this section we point out that the results presented here cannot be extended to more general piecewise-constant forcing terms.
\\ \\
\noindent {\bf Acknowledgements.} M. Zamora gratefully acknowledge support
from FONDECYT, project no. 11140203. J. Godoy was supported by a CONICYT fellowship(Chile)  in the  Program Doctorado en
Matem\'atica Aplicada, Universidad Del B\'io-B\'io no. 21161131.

\end{document}